\documentclass{article}
\usepackage{amsmath,amssymb,latexsym,misccorr,hyperref}
\author{O.\ S.\ Ogneva}
\title{Detailed proof of a theorem on coincidence of homological dimensions of Fr\'
echet algebras of smooth functions on a manifold with the
dimension of the manifold}
\date{December 14, 2012}

\newtheorem{thm}{Theorem}
\newtheorem{prop}{Propositon}
\newtheorem{cor}{Corollary}

\newcommand{\ds}{\mathop{\mathrm{ds}}}
\newcommand{\dg}{\mathop\mathrm{dg}}
\newcommand{\db}{\mathop\mathrm{db}}
\newcommand{\M}{\mathcal M}
\newcommand{\R}{\mathrm\mathbf{R}}
\newcommand{\U}{\mathcal U}
\newcommand{\C}{\mathrm\mathbf{C}}
\newcommand{\cc}{\check{C}}
\newcommand{\Co}{\mathrm\mathbf{C}\sb 0}
\newcommand{\A}{C\sp\infty (\M)}
\newcommand{\eps}{\varepsilon}
\newcommand{\al}{\alpha}
\newcommand{\ga}{\gamma}
\newcommand{\exta}{\sb{A}\mathop{\mathrm{Ext}}}
\newcommand{\ext}{\mathop{\mathrm{Ext}}}
\newcommand{\Kos}{\mathop\mathrm{Kos}}
\newcommand{\dha}{\sb A \mathop\mathrm{dh}}
\newcommand{\dsa}{\mathrm{\mathop{ds}}\, A}
\newcommand{\dga}{\mathrm{\mathop{dg}}\, A}
\newcommand{\dba}{\mathrm{\mathop{db}}\, A}
\newcommand{\prot}{\widehat\otimes}
\newcommand{\CM}{C\sp\infty(\M)}
\newcommand{\Ci}{C\sp\infty}
\newcommand{\CU}{C\sp\infty(U)}
\newcommand{\CUo}{C\sp\infty(U\sb 0)}
\newcommand{\om}{\omega}
\newcommand{\la}{\lambda}
\newcommand{\vtr}{\vartriangleright}
\newcommand{\vtl}{\vartriangleleft}
\newcommand{\di}{\mathop\mathrm{dist}}
\newcommand{\ot}{\otimes}
\renewcommand{\phi}{\varphi}

\renewcommand{\Im}{\mathop\mathrm{Im}}

\begin{document}
\maketitle
 Given work is devoted to the proof of the following
assertion.

\begin{thm}\label{t1}
 For the topological algebra $C\sp{\infty}(\M)$ of smooth
 functions on a smooth $m$-dimensional real manifold
$\M$ the small global dimension $(\ds C\sp\infty (\M))$, the
global homological dimension $(\dg C\sp\infty (\M))$ and the
bidimension $(\db C\sp\infty(\M))$ are equal to $m$ (all
dimensions are understood in the sense of the homology of
topological (locally convex) algebras~\cite{helem}).
\end{thm}

  The \textit{bidimension} of the topological algebra $C\sp\infty (U)$,
  where $U$ is an open set in $\R \sp m$, was computed by Taylor
  in the paper~\cite{tay}. The proof of this fact, like the
  proof of purely algebraic Hilbert's \textit{Sysygy} Theorem~\cite[Chapter VII, $\S 7$]{mac}
   is substantially based upon the possibility of constructing a special \textit{free resolution} of length $m$,
   the so-called \textit{Koszul resolution}. However in
   the general case of the topological algebra of smooth functions
   on an arbitrary smooth real manifold $C\sp\infty(\M)$ there is
   no natural system of commuting operators (of the type
   of the multiplication operators by an independent variable in the case of
   $\M\subset \R\sp m)$, which would let us construct free
   Koszul resolutions. Nevertheless we will show, that the modules
over the algebra $\A$ always have projective (generally  speaking,
non-free) resolutions of length $m$, however they are more
complicated than the Koszul resolutions. In the proof the
projectivity of some natural class of modules (see Section~$2$) is
employed essentially. These modules are used for the construction
of Koszul resolutions and of more complicated resolutions,
received then by means of the smooth \v{C}ech cochain complex. Let
us emphasise that the main result is related to the topological
(`locally convex') homology and its proof uses the specific
machinery of this theory.

 1. Recall some basic notions and facts, which we will use.
  Let $L$ be a locally convex space, let  $(T\sb 1, T\sb 2,\cdots,T\sb
  m)$ be a system of commuting continuous operators on $L$.
  The complex
  \[
  0\to L\otimes E\sb m\xrightarrow{d\sb m} L\otimes E\sb {m-1}
  \to\cdots\to L\otimes E\sb 1 \xrightarrow{d\sb 1} L\to 0,
  \]
  where $E\sb i = \bigwedge\sp i \C\sp m,$
  \begin{equation}\label{kos}
   d\sb i(x\otimes e\sb{j\sb 1} \wedge \cdots \wedge e\sb{j\sb i})=
  \sum\limits\sb{k=1}\sp i (-1)\sp{k-1}T\sb {j\sb k} x\otimes
  e\sb{j\sb 1}\wedge \cdots \wedge e\sb{j\sb {k-1}}\wedge
  e\sb{j\sb{k+1}}\wedge\cdots\wedge e\sb{j\sb i},
  \end{equation}
is called the \textit{Koszul complex} of the pair $(L; (T\sb 1,
\cdots, T\sb m))$ (we denote this complex by $\Kos(L;(T\sb 1,
\cdots, T\sb m))$, and the complex
\begin{equation}\label{aug}
0\to L\otimes E\sb m\xrightarrow{d\sb m} L\otimes E\sb {m-1}
  \to\cdots\to L\otimes E\sb 1 \xrightarrow{d\sb 1} L\xrightarrow{\eps} E/D\sb m\to 0,
 \end{equation}
where $D\sb m= \Im T\sb 1 + \cdots+ \Im T\sb m,$ $\eps$ is the
natural projection, is called the \textit{ augmented Koszul
complex} of the pair $(L; (T\sb 1, \cdots, T\sb m))$.

For the pair $(L; (T\sb 1, \cdots, T\sb m))$ we consider the
spaces $D\sb i= \Im T\sb 1+\cdots + \Im T\sb i$,  $i=1,\cdots, m$,
$D\sb 0=\{0\}$. Then the operator $T\sb {i+1}$, which leaves the
space $D\sb i$ invariant, induces the operator $\overline{T}\sb
{i+1}$ on $L/D\sb i$. It is known that
(see~\cite[Theorem~$V.1.3$]{helem}, \cite[Proposition~$4.1$]{tay})
an augmented Koszul complex~(\ref{aug}) has a contracting homotopy
in the category of locally convex spaces provided the operator
$\overline{T}\sb i$ has a left inverse continuous operator for all
$i=1,\cdots, m$.

 The definitions of topological (locally convex) algebras (in
 particular, Fr\'echet algebras), of modules over them, and of their
 homological characteristics are given in~\cite{helem},
 \cite{cotay}. For  a locally convex algebra $A$ we denote by $\dha X$ the \textit{homological
 dimension} of a left $A$-module $X$ (that is the minimal length
 of the projective resolution of $X$ ); by $\dsa$, $\dga$, $\dba$
 we denote, respectively, the \textit{left small homological dimension}, the
  \textit{left global dimension}, and the \textit{bidimension} of the algebra $A$.
These values are defined as follows: $\dsa=\sup\{\dha X : X \mbox{
is a left $A$-module, } \dim X < \infty \}$, $\dga=\sup\{\dha X :
X \mbox{ is a left $A$-module}\}$, $\dba = \sb{A\sp
e}\mathop\mathrm{dh} A$, where $A\sp e =A\prot A\sp{op}$ is the
enveloping algebra of the algebra $A$ (here $\prot$ is the
complete projective tensor product, and $A\sp {op}$ is the algebra
with the opposite multiplication).

  We will consider the commutative Fr{\'e}chet algebra
  $\CM$ of  smooth functions on a smooth $m$-dimensional real
  manifold $\M$. The topology on the space $\CM$ is given by the
  system of seminorms
  \[
  \|f\|\sb{(W\sb i, \om\sb i), K, n\sb 1, \cdots, n\sb
  m}=\max\limits\sb{\om\sb i(K)}  \left| \frac {\partial\sp{n\sb 1
  +\cdots+n\sb m}}{\partial x\sb 1\sp{n\sb 1}\cdots\partial x\sb m\sp{n\sb
  m}} f(\om\sb i\sp{-1}(x))\right|,
  \]
where $(W\sb i, \om\sb i)$, $i= 1, 2,\cdots$ are \textit{charts},
i.e., open sets $W\sb i\subset \M$ along with the fixed
homeomorphisms $\om\sb i$ onto open subsets in $\R\sp m$,
$K\subset W\sb i$ is a compact set, $n\sb 1,\cdots, n\sb{m}$ are
non-negative integers, $ x=(x\sb 1,\cdots, x\sb m)=\om\sb i (s)\in
\R\sp m$, $s\in W\sb i$. We recall that by Grothendieck's
Theorem~\cite[Chapitre~II, \S 3]{groth} for any smooth manifolds
$\M\sb 1$ and $\M\sb 2$ the spaces $C\sp\infty (\M\sb 1,
C\sp\infty (\M\sb 2))= C\sp\infty(\M\sb 1\times \M\sb 2)$ and
$\Ci(\M\sb 1)\prot \, \Ci(\M\sb 2)$ are topologically isomorphic.
Moreover it is evident that the isomorphism $\CM\sp e=\CM\prot \,
\CM \simeq C\sp\infty (\M\times \M)$ is an isomorphism of
algebras. 

Since (see~\cite{helem})  for arbitrary locally convex
algebra
\[
\dsa\le\dga\le \dba,
\]
in order to prove Theorem~\ref{t1} it is sufficient to find a
finite-dimensional $\CM$-module with the homological dimension not
less than $m$ (Section~$3$ is devoted to this) and to establish
that $\sb{C\sp\infty(\M)\sp e} \mathrm{\mathop{dh}} \,  C(\M)\le
m$ (see Section~$5$).

2. For an open set $U\subset \M$ the space $\CU$  of smooth
functions is a Fr{\'e}chet module over the algebra $\CM$ with
respect to the pointwise outer multiplication $f\cdot
g(s)=f(s)g(s),$ where $f\in \CM$, $g\in \CU$, $s\in U$.

\begin{thm}\label{t2}
 Let $U$ be an open set in $\M$, which is contained in a
chart. Then $\CU$ is a projective $\CM$-module.
\end{thm}

$\vtl$ It is sufficient to prove that $\CU$ is a retract of the
free $\CM$-module $\CM\prot \, \CU$,  i.e., that the canonical
projection
\[\pi\sb U : C\sp \infty (\M\times U) \to \CU,
\]
 $\pi\sb U(f) =f(s,s),$ $ s\in U$ has a left inverse morphism
$\rho$ \cite[Theorem~III.1.30]{helem}.

  Let $U\subset W,$ $(W,\om)$ be a chart. We consider the positive continuous
  function
  \[\psi(x)= \min (1, \di (x,\partial \om (U)),\]
   where  $dist$ is the Euclidean distance from the point $x\in \om(U)$ to
  the boundary.

  We take an arbitrary smooth function $\phi(x),$ $x\in \om(U)$
  such that \newline $0<\phi(x)<\psi(x)$ and, using this function, we
  define on $\om(U)\times\om(U)$
  the smooth function $\theta$ by
\[
\theta(x,y)=
\begin{cases}
\exp \frac{|x-y|\sp 2}{|x-y|\sp 2 - \phi(y)\sp 2}, &\text{if
$|x-y|\le \phi(y),$}\\0, &\text{if $|x-y|>\phi(y).$}
\end{cases}
\]

 For $(s,t)\in \M\times U$ we set
\[
F(s,t)=
\begin{cases}
\theta(\om(s),\om(t)),&\text{if $s\in U,$}\\
0, &\text{if $s\not \in  U.$ }
\end{cases}
\]
 Since it is evident that the support of $F(s,t)$ belongs to
 $U\times U$ then \newline
 $F(s,t)\in C\sp\infty (\M\times U)$, moreover
 $F(s,s)=1$. Then we can define  the required map $\rho:\,\CU \to C\sp \infty (\M\times U)$  by
 \[
 (\rho f)(s,t) =
 \begin{cases}
f(s)F(s,t), &\text{if $s\in U,$}\\
0, &\text{if $s\not\in U.$}
\end{cases}
\]
 It is clear that $\rho$ is well-defined, it is a morphism, and
 $\pi\sb{U} \rho=1\sb{C\sp\infty(U)}.$

 Thus the projectivity of the module $\CU$ is proved. $\vtr$

\medskip

 3. For any open set $U\subset \M$, lying entirely in a chart,
 e.g., $(W,\om)$, we define the multiplication operators by the `$k$-th
 coordinate function' $T\sb k\sp m(U): \CU \to \CU.$
Namely for $\om(t)=(\om\sp 1 (t),\cdots,\om\sp m(t))\in \R\sp m$,
$t\in U$ we set $T\sb k\sp m(U)f=\om\sp k (t) f(t),$ $k=1,\cdots,
m.$ It is evident that $T\sb k\sp m(U)$ are $\CM$-module
morphisms.

  We take an arbitrary point $s\sb 0\in U$ and an open  set $U\sb
  0$, containing this point and belonging to a chart.
  Without loss of generality one can suppose that $U\sb 0$
  is homeomorphic to $\R\sp m$ $(\om(U\sb 0)=\R\sp m)$ and $\om(s\sb 0)=0$.

 We denote by $\Co$ the one-dimensional $\CM$-module $\C$ with
 the outer multiplication $f\cdot \la= f(s\sb 0)\la,$ where $\la
 \in \C$, $f\in \CM.$

 \begin{thm}\label{t3}
 The complex
 \begin{equation}\label{f3}
 \Kos(\CUo; (T\sb 1\sp m(U\sb 0), \cdots, T\sb m\sp m(U\sb
 0))\xrightarrow {\pi\sb 0} \Co,
\end{equation}
over $\Co$, where $\pi\sb 0 : \CUo\to \Co$, $\pi\sb 0(f)= f(s\sb
0),$ is a projective resolution of the $\CM$-module $\Co$.
\end{thm}

$\vtl$ Since all modules $\CUo\otimes E\sb i$ are projective as
projective summands of projective modules and the maps $d\sb i$,
$\pi\sb 0$ are morphisms of $\CM$-modules it is sufficient to
prove that the complex~(\ref{f3}) is admissible, i.e., it has a
contracting homotopy in the category of Fr\'echet spaces.

It is clear that the condition $\om(U\sb 0)= \R \sp m$ implies the
isomorphism of the Fr\'echet spaces of the complex in question and
that of the complex

\begin{equation}\label{f4}
\begin{split}
 0\to \Ci(\R\sp m)\otimes E\sb
m &\xrightarrow{d\sb m} \Ci (\R\sp m)\otimes  E\sb {m-1}
\to\cdots\\
 &\to \Ci(\R\sp m)\otimes E\sb 1\xrightarrow{d\sb 1} \Ci(\R\sp
m)\xrightarrow{\pi\sb 0'} \Co\to 0
\end{split}
\end{equation}

\noindent in the category of complexes, that is the isomorphism
between (\ref{f3}) and  the complex

\[
\Kos(\Ci(\R\sp m); (T\sb 1, \cdots, T\sb m))\xrightarrow{\pi\sb
0'} \Co,
\]

\noindent where $T\sb k(f)=x\sb k f(x),$ $x=(x\sb 1,\cdots,x\sb
m)\in \R\sp m,$ $\pi\sb 0 ': \Ci (\R\sp m) \to \Co,$ $\pi\sb 0
'(f)=f(0).$

  For the complex~(\ref{f4}) the space $D\sb k = \mathrm{\mathop{Im}}\, T\sb 1+ \cdots
  +  \mathrm{\mathop{Im}}\, T\sb k$ coincides with the space $A\sb k= \{f\in \Ci(\R\sp
  m): f(0,\cdots,0,x\sb {k+1},\cdots,x\sb m)=0\}$: the inclusion
  $D\sb k \subset A\sb k$ is evident, and from Hadamard's Lemma it
   follows that functions $f\in A\sb k$ are representable in the
   form $f(x)= \sum\limits\sb{i=1}\sp k x\sb i f\sb i (x),$
   where $f\sb i\in \Ci (\R\sp m)$, i.e., they belong to $D\sb
   k$.  Therefore $\Ci(\R\sp m)/D\sb k\simeq \Ci(\R\sp{m-k})$
(in particular, $\Ci(\R\sp m)/D\sb m\simeq\Co$) and the natural
projection $\eps : \Ci(\R\sp m)\to \Ci(\R\sp m)/D\sb m$ is
identified with the morphism $\pi\sb 0'$. Thus we have showed that
up to isomorphism the complex (\ref{f4}) is the augmented Koszul
complex of the pair $(\Ci (\R\sp m);(T\sb 1, \cdots,T\sb m)),$
besides, the established isomorphisms let to consider the operator
$\overline{T}\sb {k+1}$ to be acting on the space $\Ci(\R\sp
{m-k})$:
\[
\overline{T}\sb{k+1}(g(x\sb{k+1},\cdots, x\sb m))= x\sb
{k+1}g(x\sb {k+1},\cdots,x\sb m).
\]
Then the continuous operator
\[S\sb{k+1}:
\Ci(\R\sp{m-k})\to \Ci (\R\sp{m-k}),
\]
\[
S\sb k(g(x\sb{k+1},\cdots, x\sb m))=(g(x\sb{k+1},\cdots, x\sb
m)-g(0,x\sb{k+2},\cdots, x\sb m))/x\sb{k+1}
\]
is a left inverse to $\overline{T}\sb{k+1}$. Thus the  sufficient
conditions from Section~{1} for the existence of a contracting
homotopy for the augmented Koszul complex~(\ref{f4})  in the
category of Fr\'echet spaces  are fulfilled. Therefore the
complex~(\ref{f4}) and the complex~(\ref{f3}), isomorphic to it,
are admissible. $\vtr$

\begin{prop}\label{p1}
$\sb {\CM} \mathop{\mathrm {dh}}\Co =m.$
\end{prop}
$\vtl$
 In view of the fact that \cite[Theorem~$III.5.4$]{helem} for an arbitrary
 locally convex algebra $A$  and for an $A$-module $X$
\begin{align*}
 \dha X = \sup \{k: {\exta}\sp{k+n}(X,Y)=0  \text{ for any
$n>0$ and }\\
\text { there is a module $Y$, such that } {\exta}\sp k(X,Y)\not=0
\}
\end{align*}
it is sufficient to prove that ${\sb{\CM}\mathop{\mathrm{Ext}}}\sp
m (\Co,\Co)\not=0$. For this  we use the projective resolution
(\ref{f3}) of the module $\Co$. Up to an isomorphism the
 left end of the resolution has the form
 \[
 0\to \CUo \xrightarrow{d\sb m} \underbrace{\CUo\oplus\cdots\oplus\CUo}_m\to\cdots,
 \]
where $d\sb m(f)=(\om\sp 1(t)f(t),\cdots,\om\sp m(t)f(t)),$ $t\in
U\sb 0$. Since for any function $f\in \CUo$ there exists a
sequence $\{f\sb n\}$, $f_n \in \CM$, such that $f=\lim\limits\sb{n\to \infty}
f\sb n|\sb {U\sb 0}$ in $\CUo$, then for any $\CM$-module morphism
$\al: \CUo\to \Co$ we have
\begin{align*}
 \al(f)=\al
(\lim\sb{n\to\infty}f\sb n|\sb{U\sb
0})&=\lim\sb{n\to\infty}\al(f\sb n|\sb{U\sb 0}) \\
=\lim\sb{n\to\infty}(f\sb n\cdot \al(1))&= \lim\sb{n\to\infty}f\sb
n(s\sb 0)\al(1)=f(s\sb 0)\al(1).
\end{align*}
Consequently,  for any morphism $\ga:
\underbrace{\CUo\oplus\cdots\oplus\CUo}_m\to \Co$ we have
\begin{align*}
\ga d\sb m(f)&=\ga(\om\sp 1(t)f(t),\cdots, \om\sp
m(t)f(t))\\
&=\ga(\om\sp 1(t)f(t))+\cdots+\ga(\om\sp m (t))f(t)=0
\end{align*}
and then $\sb{\CM} {\mathop{\mathrm{Ext}}}
(\Co,\Co)=\sb{\CM}\mathop\mathrm{hom}(\CUo,\Co).$

It remains to observe that the later space is not zero because it
contains the non-trivial morphism $\pi\sb0$, defined by $\pi\sb
0(f)=f(s\sb 0)$.$\vtr$

\medskip

 From  Proposition~\ref{p1} we immediately get a lower
 estimate for homological dimensions.

\medskip

\begin{cor}
$\mathop\mathrm{ds}\CM\ge m.$
\end{cor}

\medskip

4. The aim of the further exposition is to get an upper estimate
for homological dimensions. At the beginning we make it `locally',
namely for the $\CM\sp e$-module $\CU$ with the outer
multiplication $f\cdot g= f(s,s)g(s),$ $f\in \CM\sp e,$ $g\in
\CU,$ $s\in U,$ we establish that $\sb{\CM\sp
e}\mathrm{\mathop{dh}}\, \CU\le m.$

 As before we suppose that $U$ is contained in a chart.

\begin{thm}\label{t4}
 The complex
 \begin{equation}\label{f5}
  \begin{split}
 \Kos\, (\Ci(U\times U);(T\sb 1\sp{2m}(U\times
 U)-T\sb{m+1}\sp{2m}(U\times U),\cdots,\\
 T\sb{m}\sp{2m}(U\times U)-T\sb{2m}\sp{2m}(U\times U)))
 \xrightarrow{\pi}\CU,
  \end{split}
\end{equation}
where $\pi(f)=f(s,s),$ $s\in U$  is a projective resolution of
$\CM\sp e$-module $\CU$.
\end{thm}

$\vtl$
 By Theorem~\ref{t2} the $\CM\sp e$-modules $\Ci(U\times
U)\otimes E\sb i$ are projective and consequently, as in
Theorem~\ref{t3} it is sufficient to establish the admissibility
of the complex~(\ref{f5}).

 But the complex~(\ref{f5}) is isomorphic to the complex

\begin{equation}\label{f6}
\begin{split}
0\to \Ci (\om(U)\times\om(U))\otimes E\sb m\xrightarrow{d\sb
m}\Ci(\om(U)\times \om(U))\otimes E\sb{m-1}\to \cdots \\
\to \Ci(\om(U)\times\om(U))\xrightarrow{\pi'}\Ci(\om(U))\to 0,
\end{split}
\end{equation}

$(\om(U)\subset \R\sp m),$
 which is known to be admissible \cite[Proposition~$4.4$]{tay}.
$\vtr$

\medskip

\begin{cor}
$\sb{\CM\sp e} \mathop{\mathrm{dh}}\CU\le m.$
\end{cor}

\medskip

We take an arbitrary set $\{U\sb i\},$ $i=1,2,\cdots$ of open
subsets of $\M$, each of which lies in a chart. Then the space
$\prod\limits\sb{i=1}\sp{\infty} \Ci(U\sb i)$ being a countable
Cartesian product of Fr\'echet spaces and at the same time being
$\CM\sp e$-modules is a Fr\'echet  $\CM\sp e$-module \cite{helem}
with respect to the outer multiplication, defined componentwise.

\begin{prop}\label{p2}
$\sb{\CM\sp e}\mathrm{\mathop{dh}} \prod\limits\sb{i=1}\sp\infty
\Ci (U\sb i)\le m.$
\end{prop}

$\vtl$ We consider the Cartesian product of the
complexes~(\ref{f5}) for \newline  $\CM\sp e$-modules   $\Ci(U\sb
i)$:

\begin{equation}\label{f7}
\begin{split}
0\to \prod\limits \sb{i=1}\sp \infty \Ci(U\sb i\times U\sb
i)&\otimes E\sb m\xrightarrow{d\sb m}
\prod\limits\sb{i=1}\sp\infty \Ci(U\sb i\times U\sb
i)\otimes E\sb{m-1}\xrightarrow{d\sb{m-1}} \cdots \\
 \to
&\prod\limits\sb{i=1}\sp\infty \Ci(U\sb i\times U\sb
i)\xrightarrow{\pi} \prod\limits \sb{i=1}\sp \infty \Ci(U\sb i)\to
0.
\end{split}
\end{equation}

 We prove that the complex~(\ref{f7}) is a projective resolution
 for $\prod\limits\sb{i=1}\sp \infty \Ci(U\sb i)$.

Actually, it is evident that the maps $d\sb i$, $\pi$ are $\CM\sp
e$-module morphisms. The contracting homotopy maps for the
complex~(\ref{f7}) are obtained as Cartesian products of
respective contracting homotopy maps for $\CM\sp e$-module
complexes for $\Ci (U\sb i)$, therefore the admissibility of the
complex~(\ref{f7}) follows from the admissibility of the
complex~(\ref{f5}).

   It remains to show that the modules
   $\prod\limits\sb{i=1}\sp\infty \Ci (U\sb i\times U\sb i)$ are
   projective. As in Theorem~\ref{t2} for the canonical
   morphism
   \[
   \pi\sb \infty: \Ci (\M\times \M) \prot \prod\limits
   \sb{i=1}\sp\infty \Ci(U\sb i\times U\sb i) \to
   \prod\limits\sb{i=1}\sp \infty
   \Ci(U\sb i\times U\sb i),
   \]
\[
\pi\sb\infty(f\otimes\{g\sb i\}\sb {i=1}\sp \infty)=\{\pi\sb
i(f\ot g\sb i)\}\sb{i=1}\sp\infty
\]

(here $\pi\sb i$ are canonical $\CM\sp e$-module morphism for the
module ${\Ci (U\sb i\times U\sb i)}$) we suggest a  morhism
$\rho\sb\infty$ such that $\pi\sb \infty \rho\sb\infty=1$. For the
constructing the morphism $\rho$ we use the Fr\'echet space
isomorphism~\cite[Chapter~$1$, $\S1$, Proposition~$1$]{groth}
\[
\Ci (\M\times \M)\prot \prod\limits\sb{i=1}\sp\infty \Ci(U\sb
i\times U\sb i)\simeq \prod\limits\sb{i=1}\sp\infty \Ci(\M\times
\M)\prot Ci(U\sb i\times U\sb i)
\]
 and we define the morphism $\rho\sb \infty$ in the following
 manner
\[
\rho\sb\infty:\prod\limits\sb{i=1}\sp \infty \Ci(U\sb i\times U\sb
i)\to \prod\limits \sb{i=1}\sp\infty \Ci (\M\times \M)\prot
\Ci(U\sb i\times U\sb i),
\]
\[
\rho\sb\infty (\{g\sb i\}\sb{i=1}\sp\infty)=\{\rho\sb i g\sb
i\}\sb{i=1}\sp\infty
\]
($\rho\sb i$ is a left inverse $\CM\sp e$-module morphism to the
morphism $\pi\sb i$, which exists because of projectivity of the
module $\Ci(U\sb i\times U\sb i$).

 Thus the length of the admissible resolution~(\ref{f7}) for the
 $\CM\sp e$-module $\prod\limits\sb{i=1}\sp\infty \Ci(U\sb i)$ does not
 exceed $m$ and consequently
 $\sb{\CM\sp e}\mathrm{\mathop{dh}} \prod\limits\sb{i=1}\sp\infty
 \Ci(U\sb i)\le m.\vtr$

 \medskip

5.\  The final phase of the proof is to establish the estimate
$\sb{\CM\sp e}\mathrm{\mathop{dh}} \CM\le m$ by `pasting together'
the local upper estimates by the means of the \v{C}ech complex.

\begin{thm}\label{t5}
$\sb{\CM\sp e}\mathrm{\mathop{dh}} \CM\le m$
\end{thm}
$\vtl$ We take a covering of the manifold $\M$ with countable set
of charts $\{(W\sb i, \om\sb i)\}$. The covering dimension of
$m$-dimensional smooth real  manifold does not exceed $m$ (see
\cite[Theorem~$2.15$]{man}), therefore  one can inscribe a
countable covering $\U=\{U\sb i\}$ of the multiplicity less at
most $m$ in the covering $\{W\sb i\}$. In other words any point of
the manifold is contained in the $m+1$ sets of the system $\U$ at
most. We construct the augmented smooth cochains complex
corresponding to the covering $\U$:
\begin{equation}\label{f8}
\begin{split}
 0\to
\CM\xrightarrow{\eta} &\cc\sp 0(\U,\CM) \xrightarrow{\partial\sp
0}\\
&\cc\sp 1(\U,\CM)\xrightarrow{\partial\sp 1}\cdots
\xrightarrow{\partial\sp {m-1}} \cc\sp m(\U,\CM)\to 0,
\end{split}
\end{equation}
where $\cc\sp i(\U,\CM)=\prod\limits\sb{\sigma\sb i}
\Ci(|\sigma\sb i|),$

(the product is taken by the sets $\sigma\sb i=(U\sb{j\sb
0},\cdots, U\sb{j\sb i}),$ ${j\sb 0}<\cdots< {j\sb i}$ such that
$|\sigma \sb i|=U\sb{j\sb 0}\cap\cdots \cap U\sb{j\sb i}\not
=\emptyset$)

\[
(\partial\sp{i-1} f)(\sigma\sb i)=\sum\sb{k=0}\sp i (-1)\sp k
f(U\sb{j\sb
0},\cdots,U\sb{j\sb{k-1}},U\sb{j\sb{k+1}},\cdots,U\sb{j\sb
i})|\sb{|\sigma\sb i|},
\]

\[\eta f=\{f|\sb {U\sb i}\}\sb {i=1}\sp\infty.\]

 We note that $\cc\sp i (\U,\CM)=0$ for $i> m$ because $|\sigma\sb
 i|=\emptyset.$

 The complex~(\ref{f8}) is exact~\cite{gan}. Moreover, one can
 remark that in the proof of this fact the necessary contracting homotopy in the category of
 Fr\'echet spaces  is constructed
 \cite[Lemma~$VI.D.3$]{gan}.
Therefore the complex~(\ref{f8}) is admissible. We represent the
complex in the form of Yoneda product of $m+1$ short exact
admissible sequences of $\CM\sp e$-modules and morphisms:

\begin{equation}\label{f9}
0\to \mathrm{\mathop{Ker}}\, \partial\sp {m-1}\to
\cc\sp{m-1}(\U,\CM)\xrightarrow{\partial\sp{m-1}}\cc\sp
m(\U,\CM)\to 0,
\end{equation}

\begin{equation}\label{f10}
0\to \mathrm{\mathop{Ker}}\, \partial\sp {i}\to
\cc\sp{i}(\U,\CM)\xrightarrow{\partial\sp i}
\mathrm{\mathop{Ker}}\,
\partial\sp{i-1}\to 0,
\end{equation}
\[
i=1,\cdots,m-2,
\]
\begin{equation}\label{f11}
0\to \CM \xrightarrow{\eta}\cc\sp 0
(\U,\CM)\xrightarrow{\partial\sp 0} \mathrm{\mathop{Ker}}\,
\partial\sp 1 \to 0.
\end{equation}

 For an arbitrary fixed $\CM\sp e$-module $Y$ each of the
 sequences (\ref{f9})--(\ref{f11}) defines a long exact
 sequence~\cite[Theorm~$III.4.4$]{helem}.
 For instance, for the sequence~(\ref{f9}) the long exact sequence
 has the form:
 \begin{align*}
\cdots\to\sb{\CM\sp e}{\ext}\sp k(\cc\sp{m-1}, Y)\to \sb{\CM\sp
e}{\ext}\sp k (\mathop{\mathrm{Ker}}\, \partial\sp{m-1}, Y)\\ \to
\sb{\CM\sp e}{\ext}\sp {k+1}(\cc\sp m, Y)\to \cdots.
\end{align*}

Since $\sb{\CM\sp e} \mathop\mathrm{dh} \cc\sp{m-1}\le m$ and
$\sb{\CM\sp e} \mathop\mathrm{dh} \cc\sp{m}\le m$, then
\[{\sb{\CM\sp e}{\ext}\sp k (\mathrm{\mathop{Ker}}\,
\partial\sp{m-1}, Y)=0}
\]
 when $k>m$. Applying further the exact sequence for the functor
$\sb{\CM\sp e}{\ext}(\cdot, Y)$ to the sequence (\ref{f10}) when
$i=m-2,\cdots,1,$ and then to the sequence~(\ref{f11}), and,
taking into account Proposition~\ref{p2}, we get that $\sb{\CM\sp
e}{\ext}\sp k (\mathrm{\mathop{Ker}}\,
\partial\sp i, Y)=0$ when $k>m$, and finally, $\sb{\CM\sp e}
{\ext}\sp k (\CM, Y)=0$ when $k>m$. Consequently, $\sb{\CM\sp
e}\mathrm{\mathop{dh}} \CM \le m.$ $\vtr$

\medskip

  Joining the just obtained inequality and
  Corollary from Theorem~\ref{t3}
we get that
\[
\ds \CM=\dg \CM=\db \CM=m.
\]
 The proof of Theorem~\ref{t1} is finished.


\begin{thebibliography}{XX}
\bibitem{helem}
A.\ Ya.\ Helemskii, The homology of Banach and topological
algebras. Kluwer Academic Publishers, Dordrecht (1989).

\bibitem{tay}
J.\ L.\ Taylor, A general framework for a multioperator functional
calculus. Advances in Math., 9(1972), 183--252.

\bibitem{mac}
S.\ MacLane, Homology. Springer-Verlag, Berlin, 1963.

\bibitem{cotay}
J.\ L.\ Taylor, Homology and cohomology for topological algebras.
Advances in Math., 9(1973), 137--182.

\bibitem{groth}
A.\ Grothendieck, Produits tensoriels topologiques et espaces
nucl\'{e}aires. Mem.\ Amer.\ Math.\ Soc., 16(1955).

\bibitem{man}
J.\ R.\ Munkres, Elementary differential topology. Revised ed.
Ann.\ Math.\ Studies, 54, Princeton Univetsity Press, 1966.

\bibitem{gan}
R.\ C.\ Gunning, H.\ Rossy, Anaytic functions of several complex
variables. Englewoods Cliffs, N.\ J.\
, 1965.







\end{thebibliography}
\end{document}